\documentclass[a4paper,12pt]{article}
\usepackage{amssymb,amsmath}
\usepackage[dvips]{color}
\usepackage[dvips]{graphics}
\newtheorem{theorem}{Theorem}[section]
\newtheorem{lemma}{Lemma}[section]

\newtheorem{example}{Example}

\newcommand{\eps}{\varepsilon_M}

\begin{document}
\title{Numerical stability of iterative refinement with a
relaxation for linear systems}
\author{Alicja Smoktunowicz, Jakub Kierzkowski and Iwona
Wr\'obel \thanks{ Faculty of Mathematics and Information
Science,
Warsaw University of Technology, Koszykowa 75, 00-662 Warsaw,
Poland,
e-mails: smok@mini.pw.edu.pl, J.Kierzkowski@mini.pw.edu.pl,
wrubelki@wp.pl }}
  \maketitle

\begin{abstract}
Stability analysis of Wilkinson's iterative refinement with a
relaxation IR($\omega$) for solving linear systems 
is given.
It extends existing results for $\omega=1$, i.e., for
Wilkinson's iterative refinement. We assume that all
computations are performed in fixed (working) precision
arithmetic.
Numerical tests were done in \textsl{MATLAB} to illustrate our
theoretical results.
A particular emphasis is given on convergence of iterative
refinement with a relaxation. Our tests confirm that the choice
$\omega=1$ is the best choice from the point of numerical
stability.
\end{abstract}

\medskip


\medskip

\noindent {\bf AMS2010 MSC:} 65F05,  65F10, 15A12

\section{Introduction}
\label{intro}
We consider the system $Ax=b$, where $A \in \mathbb R^{n \times
n}$ is nonsingular and $b \in \mathbb R^{n}$. Iterative
refinement techniques for linear systems of equations are very
useful in practice and the literature on this subject is very
rich  (see \cite{arioli}--\cite{bjorck3}, \cite{higham1}-- \cite{Zlatev}).

The idea of relaxing the iterative refinement step is the
following. We require a basic linear equation solver $S$ for
$Ax=b$ which uses a factorization of A into simple factors
(e.g., triangular, block-triangular etc.).
Such factorization is used again in the next steps of iterative
refinement. Wilkinson's iterative refinement with a relaxation
IR($\omega$) consists of three steps.

\medskip

{\bf Algorithm  IR($\omega$) }

\medskip

Given $\omega>0$. Let $x_0$ be computed by the solver $ S$.

For $k=0,1, 2, \ldots $, the $k$th iteration consists of the
three steps:
\begin{enumerate}
\item Compute $r_k=b-A x_k$.
\item Solve $A p_k=r_k$ for $p_k$ by the basic solution solver
$S$.
\item Add the correction, $x_{k+1}=x_k + \omega \,  p_k$.
\end{enumerate}

\medskip

Clearly, $\omega=1$ corresponds to Wilkinson's iterative
refinement (see \cite{Wilkinson}).  Wu and  Wang (see \cite{Wu})
proposed this method for $\omega=\frac{h}{h+1}$, where $h>0$
(i.e., for $0<\omega<1$).
They developed the method as the numerical integration of a
dynamic system with step size $h$. 
A  preliminay error analysis of the Algorithm IR($\omega)$ was given  in \cite{Wu}  for  $0<  \omega <1$, assuming that
the extended precision is used for computing the residual
vectors $r_k$. Wu and  Wang considered only Gaussian
elimination as a solver $S$ .

The purpose of this paper is to analyze the convergence of this
method for $0<  \omega < 2$   and to show with examples that  the choice
$\omega=1$ is the best choice from the point of numerical
stability.

Notice that for arbitrary $\omega>0 $, the IR($\omega$) method
is a stationary method (in the theory) and we have $p_k=A^{-1}
r_k= x^{*}- x_k$, so $x_{k+1}- x^{*}= (1-\omega) ( x_{k}-
x^{*}), \quad k=0,1, \ldots$,
where $x^{*}$ is the exact solution to $Ax=b$. We see that the
sequence $\{x_k\}$ is convergent for arbitrary initial $x_0$ if
and only if $0 < \omega <2$.
For $\omega=1$ (Wilkinson's iterative refinement) $x_1$ will be the exact solution 
 $x^{*}$. 
It is interesting to check the influence on the relaxation
parameter $\omega$ on numerical properties of the algorithm
IR($\omega)$, assuming that all computations are performed only
in the working (fixed) precision.

Throughout the paper we use only the $2$-norm and assume that
all computations are performed in the working (fixed) precision.
We use a floating point arithmetic which satisfies the IEEE
floating point standard (see Chapter~2 in \cite{higham:96}).
For  two floating point numbers $a$ and $b$  we have
\[
f\ell(a \diamond b) = (a \diamond b) (1+ \Delta), \quad |\Delta|
\leq \eps
\]
for results in the normalized range, where $\diamond$ denotes
any of the elementary scalar
operations  $+, -, *, /$ and $\eps$ is machine precision.

In this paper we present a comparison of Wilkinson's iterative
refinement with a relaxation IR($\omega$)   from
the point of view of numerical stability.
More precisely, we say that the computed $\tilde x$ in floating
point arithmetic is a {\bf forward stable} solution to $Ax=b$ if
\begin{equation}\label{stability}
\|\tilde x - x^{*}\| \leq \mathcal O (\eps) \kappa(A) \|x^{*}\|.
\end{equation}
Throughout this paper, $\| \cdot \|$ is the matrix or vector
two--norm depending upon context, and $\kappa(A)= \|A^{-1}\|
\|A\| $ denotes the standard condition number of the matrix $A$.

A stronger property than forward stability is backward
stability. It means that the computed $\tilde x$ in floating
point arithmetic is the exact solution of a slightly perturbed
system
\begin{equation}\label{backward_stability}
(A+\Delta A) \tilde x =b, \quad \| \Delta A\| \leq \mathcal O
(\eps) \| A\|.
\end{equation}

Our analysis is similar in spirit to \cite{higham1}--\cite{jankowski} and  \cite{skeel}-\cite{smok3}.
Jankowski and Wo\'{z}niakowski  proved in  \cite{jankowski}
that an arbitrary solver $S$ which satisfies (\ref{S}),
supported by
iterative refinement, is normwise forward stable as long as $A$
is not too ill-conditioned (say, $\eps \kappa(A) <1$), and is
normwise backward stable under additional condition $q \kappa(A)
< 1$.
We extend their results for the algorithm IR($\omega)$. 

The paper is organized as follows. A proof of numerical
stability of IR($\omega)$ is given in Section 2.
In Section 3, we present some numerical experiments that
illustrate our theoretical results.

\section{Forward stability of  IR($\omega)$}
\label{forward_stab}
We require a basic linear equation solver $S$ for $Ax=b$ such
that the computed solution $\tilde x$ by $S$ satisfies
\begin{equation}\label{S}
\|\tilde x -x^{*}\| \leq q \,\, \| x^{*}|\|, \quad q \leq 0.1.
\end{equation}
 We make a standard assumption that the matrix-vector
multiplication is backward stable. Then the computed residual
vector $\tilde r =f\ell(b-A \tilde x)$ satisfies
\begin{equation}\label{res}
\tilde r =b-A \tilde x + \Delta r, \quad \| \Delta r \| \leq
L(n) \eps ( ||b||+ ||A|| ||\tilde x||),
\end{equation}
where $L(n)$ is a modestly growing function on $n$. 

We start with the following lemma. 
\begin{lemma}\label{lemat1}
Let IR($\omega)$ for $\omega \in (0,2)$ be applied to the
nonsingular linear system $Ax=b$ using the solver $S$ satisfying
(\ref{S})-(\ref{res}).
Let $\tilde x_k$, $\tilde r_k$ and $\tilde p_k$ denote the
computed vectors in floating point arithmetic. Assume that
\begin{equation}\label{zal1}
\eps  \leq 0.01, \quad   L(n) \eps\,  \kappa(A) \leq 0.01
\end{equation}
and
\begin{equation}\label{zal2}
|1-\omega| +  \omega q  \leq 0.6.
\end{equation}

Then for $k=0,1, \ldots $ we have 
\begin{equation}\label{teza1}
||\tilde x_k - x^{*}|| \leq q_k  ||x^{*}||, \quad  q_k \leq 0.1,\end{equation}
where
\begin{equation}\label{teza2}
q_{k+1}= (|1-\omega| + q \omega) q_k + 2.31 \omega L (n) \eps
\kappa(A)+ 1.64 \eps,
\end{equation}
with $q_0=q$.
\end{lemma}

\medskip

\begin{proof} Assume that (\ref{teza1}) holds for $k$. We prove
that it holds also for $k+1$, i.e.
$||\tilde x_{k+1} - x^{*}|| \leq q_{k+1} \,\, ||x^{*}||$, where
$q_{k+1} \leq 0.1$ and
$q_{k+1}$ satisfies (\ref{teza2}).

Under assumption (\ref{res}), the computed vectors $\tilde r_k$
satisfy
\begin{equation}\label{eqs1}
\tilde r_k= b- A \tilde x_k + \Delta r_k, \quad || \Delta r_k||
\leq \eps L(n) ( ||b||+ ||A|| ||\tilde x_k||).
\end{equation}
Under assumption (\ref{S})  we have 
\begin{equation}\label{eqs2}
\tilde p_k= p_k^{*}+ \Delta p_k, \quad p_k^{*}= A^{-1} \tilde
r_k, \quad \|\Delta p_k \| \leq q \|p_k^{*}\|.
\end{equation}
Standard error analysis shows
\begin{equation}\label{eqs3}
\tilde x_{k+1}= (I+D_k^{(1)})(\tilde x_k+ (I+D_k^{(2)} ) \omega
\tilde p_k), \quad \| D_k^{(i)}\| \leq \eps.
\end{equation}

By inductive assertion, we have $\| \tilde x_k - x^{*}\| \leq
q_k \| x^{*}\| $. Hence
$$
\|\tilde x_k \|= \| x^{*} + (\tilde x_k- x^{*}) \| \leq \|
x^{*}\| +\| \tilde x_k- x^{*}\| \leq (1+q_k) \| x^{*}\|.
$$
Similarly, from (\ref{eqs2}) it follows that $\|\tilde p_k \|
\leq (1+q) \| p_k^{*}\|$, thus
\begin{equation}\label{eqs4}
\|\tilde x_k \| \leq 1.1 \| x^{*}\|, \quad \|\tilde p_k \| \leq
1.1 \| p_k^{*}\|.
\end{equation}

From (\ref{eqs1}) and the inequality $\|b\|=\|A x_{*}\| \leq
\|A\| \| x_{*}\|$ it can be seen that
\begin{equation}\label{eqs44}
\tilde r_k= b- A \tilde x_k + \Delta r_k, \quad \| \Delta r_k \|
\leq 2.1 L(n) \eps \|A\| \| x^{*}\|.
\end{equation}

We have 
\begin{equation}\label{eqs5}
p_k^{*}= A^{-1} \tilde r_k = x^{*}- \tilde x_k + \xi_k, \quad
\xi_k= A^{-1} \Delta r_k.
\end{equation}

This together with  (\ref{eqs44}) implies  the bounds 
\begin{equation}\label{eqs6}
\| p_k^{*} \| \leq \|\tilde x_k -x^{*}\| +\|\xi_k\|, \quad \|
\xi_k \| \leq 2.1 L(n) \eps \kappa(A) \| x^{*}\|.
\end{equation}

Now our task is to bound the error $\|\tilde x_{k+1} - x^{*}\|$.
For simplicity, we define $D_k^{(3)}$ such that
$$
I+D_k^{(3)}=(I+D_k^{(1)})(I+D_k^{(2)}).
$$ 

Clearly, $\| D_k^{(3)}\| \leq 2 \eps+ \eps^2$, so from
(\ref{eqs3}) we get
\begin{equation}\label{eqs7}
\tilde x_{k+1}= (\tilde x_k+ \omega \tilde p_k) + \eta_k, \quad
\|\eta_k\| \leq \eps \| \tilde x_k\| + ( 2 \eps+ \eps^2) \omega
\| \tilde p_k\|.
\end{equation}
This together with (\ref{eqs2}) and (\ref{eqs5}) gives the
identity
\[
\tilde x_{k+1}-x^{*}= (1-\omega) (\tilde x_k - x^{*}) + \eta_k+
\omega (\xi_k + \Delta p_k).
\]
Taking norms  and using (\ref{eqs2}), we obtain
\begin{equation}\label{eqs8}
\| \tilde x_{k+1}-x^{*}\| \leq |1-\omega| \| \tilde x_k -
x^{*}\| + \|\eta_k\| + \omega \|\xi_k\| + \omega q \| p_k^{*}
\|.
\end{equation}

First we estimate  $\|\eta_k\|$. 
Since $ \|\tilde x_k -x^{*}\| \leq 0.1 \| x^{*}\|$, so by
assumption (\ref{zal1}) we obtain from (\ref{eqs6}) the bounds
\begin{equation}\label{eqs6a}
\| \xi_k \| \leq 0.021 \| x^{*}\|, \quad \| p_k^{*} \| \leq
0.121 \| x^{*}\|.
\end{equation}
From (\ref{eqs4}) and (\ref{eqs7}) we have $\|\eta_k\| \leq 1.1
\eps ( \|x^{*}\| + (2+\eps) \omega \| p_k^{*} \|)$.
Now we apply (\ref{zal1}) and (\ref{eqs6a}). Since $\omega <2$,
we see that $\|\eta_k \| \leq 1.64 \, \eps \|x^{*} \|$.
Therefore,
\[
\omega \|\xi_k\| +\|\eta_k\| \leq \omega 2.1 L(n) \eps \kappa(A)
\|x^{*} \| + 1.64 \, \eps \|x^{*} \|
\]
and by (\ref{eqs6}) we get 
\[
\omega q \| p_k^{*} \| \leq \omega q \|\tilde x_k -x^{*}\| +
\omega q 2.1 L(n) \eps \kappa(A) \| x^{*}\|.
\]
Hence, from (\ref{eqs8}) and by (\ref{zal1})-(\ref{zal2}) we
finally obtain
\[ 
\| \tilde x_{k+1}-x^{*}\| \leq (|1-\omega| + \omega q) \| \tilde
x_k - x^{*}\| + 2.31 \omega L (n) \eps \kappa(A) + 1.64 \eps
\|x^{*} \|.
\]

We conclude that $\|\tilde x_{k+1}- x^{*} \| \leq q_{k+1}
\|x^{*} \|$, with $q_{k+1}$ defined in (\ref{teza2}). It remains
to prove that $q_{k+1} \leq 0.1$.
By assumptions (\ref{zal1}) and (\ref{zal2}) and using the fact
that $q_k \leq 0.1$, we see that $q_{k+1} \leq 0.6 *0.1 +
(0.0231+ 0.0164)$, so $q_{k+1} \leq 0.1$. This completes the
proof.
\end{proof}

\medskip

\begin{theorem}\label{thm1}
Under the assumptions of Lemma \ref{lemat1} the algorithm
IR($\omega)$ is forward stable for $\omega \in (0,2)$.
There exists $k^{*}$ depending only on $n$ such that for every
$k \ge k^{*}$ the following inequality holds
\begin{equation}\label{teza3}
\|\tilde x_k - x^{*} \| \leq (11.6 L(n)+ 4.2) \eps \,\,
\kappa(A) \|x^{*} \|.
\end{equation}
\end{theorem}

\medskip

\begin{proof}
We apply the results of Lemma \ref{lemat1}. Notice that from
(\ref{teza1})-(\ref{teza2}) and by assumptions (\ref{zal1}) it
follows that
\[
q_{k+1} \leq q_k 0.6 + 2.31 \omega L (n) \eps \kappa(A) + 1.64
\eps.
\]
Since $\omega < 2$ and $1 \leq \kappa(A)$, we get 
\[
q_{k+1} \leq q_k 0.6 + (4.62 L(n)+ 1.64)  \eps  \kappa(A).
\]
From this it follows that
\[
q_{k+1} \leq (0.6)^k + \frac{4.62 L(n)+ 1.64}{1-0.6} \eps
\kappa(A).
\]
From this (\ref{teza3}) follows immediately.
\end{proof}

\medskip

\section{Numerical tests}
In this section we present numerical experiments that show the 
comparison of the IR($\omega)$ for different values of $\omega$.
All tests were performed in \textsl{MATLAB} version
8\textsl{.4.0.150421 (R2014b)}, with $\eps \approx 2.2 \cdot
10^{-16}$.

Let $x^{*}= A^{-1} b$ be the exact solution to $Ax=b$ and let
$\tilde x_k$ be the computed approximation to $x^{*}$ by
IR($\omega)$.
We produced the $n\times n$ matrix $A$ and the vector $b= A
x^{*}$, with $x^{*}= (1,1, \ldots, 1)^T$.

We report the following statistics for each iteration:
\begin{itemize}
\item forward stability  error 
\begin{equation}\label{fo}
\alpha(A,b,\tilde x_k)= \frac{\|\tilde x_k -
x^{*}\|}{\kappa(A)\, \|x^{*}\|} ,
\end{equation}
\item backward stability error 
\begin{equation}\label{ba}
\beta(A,b,\tilde x_k)= \frac{\|b- A \tilde x_k\|}{\|A\| \,
\|\tilde x_k\|},
\end{equation}
\item componentwise backward stability error 
\begin{equation}\label{co}
\gamma(A,b,\tilde x_k)= \max_i {\frac{(|b- A \tilde x_k|)_i}{
(\vert A \vert \, \vert \tilde x_k \vert)_i}}.
\end{equation}
\end{itemize} 

Note that, the componentwise stability implies the backward
stability, and backward stability implies forward stability.

\medskip

We consider the following  solvers  $S$.

\begin{itemize}
\item [] {\bf Algorithm I (GEPP).}
Gaussian elimination with partial pivoting (GEPP) for the system
$Ax=b$.

\medskip 
 
\item[] {\bf Algorithm II (BLU).}
This method uses a block $LU$ factorization $A=LU$
(\cite{Demmel: 95}):
\begin{equation}\label{block2}
A=\left(
\begin{array}{cc}
A_{11} & A_{12} \\
A_{21}  & A_{22}
\end{array}
\right)= 
\left(
\begin{array}{cc}
I  & 0 \\
L_{21}  & I
\end{array}
\right)
\quad
\left(
\begin{array}{cc}
U_{11} & U_{12} \\
0   & U_{22}
\end{array}
\right).
\end{equation}
We assume that  $A_{11}(m \times m)$ is nonsingular. Then 

\begin{enumerate}
\item $U_{11}=A_{11}$, $U_{12}=A_{12}$.
\item Solve the system $L_{21} A_{11}=A_{21}$ for $L_{21}$ (by
GEPP).
\item Compute the Schur complement $U_{22} = A_{22}-L_{21}
A_{12}$.
\end{enumerate}
\end{itemize}

Next we solve the system $L U x=b$ by solving two linear
systems, using the \textsl{MATLAB} commands
\begin{verbatim}
y=L\b;  x=U\y;
\end{verbatim}

\medskip

\begin{example} \label{WilkGEPP}
Take $A=W_n$, where $W_n$ is the famous Wilkinson's matrix of
order $n$:
\begin{equation}\label{Wilk100}
W_n=\begin{pmatrix} 
1 & 0 & 0 & \dots & 0 & 1\\
-1 & 1 & 0 & \dots & 0 & 1 \\
-1 & -1 & 1 & \dots & 0 & 1\\
\vdots & \vdots & \vdots & \ddots & \ddots &  \vdots \\
-1 & -1 & -1 & \dots & -1 &   1
\end{pmatrix}.
\end{equation}

We cite R.D. Skeel who wrote in \cite{skeel}: "Gaussian
elimination with pivoting is not always as accurate as one might
reasonably expect".
It is known (see \cite{Wilkinson}) that GEPP is considered
numerically stable unless the growth factor $\rho_n$ is large.
For Wilkinson's matrix $W_n$ we have $\rho_n =2^{n-1}$.  
Description of other types of matrices for which the growth
factor is very large is given in \cite{Foster} and
\cite{higham:96}, Section $9$. 
It is interesting that for $n=100$ the Wilkinson matrix is
perfectly well-conditioned, but GEPP produces an unstable
solution!
After one step of Wilkinson's iterative refinement (for
$\omega=1$) we get the exact solution $x^{*}= (1,1, \ldots,
1)^T$.
The situation is different for other choices of parameter
$\omega$. The results are contained in Table \ref{wilk_fo}.

\begin{table}
\caption{Values of the forward stability error (\ref{fo}) for
Algorithm I (GEPP) , where $A$ is the
$100\times 100$ Wilkinson matrix defined in (\ref{Wilk100}).
Here $\kappa(A)=44.8$. }
\label{wilk_fo}
\begin{tabular}{lllllll}
\hline $\omega  \slash k$ & %
$0.3$ & $0.5$ & $0.7$ & %
$0.9$ & $1.0$ & %
$1.2$ \\
\hline 
$0\hskip-.6ex$ & 
1.51E-02 & 1.51E-02 & 1.51E-02 & 1.51E-02 & 1.51E-02 & 1.51E-02
\\
$1\hskip-.6ex$ &
1.05E-02 &	7.56E-03 &	4.54E-03 &	1.51E-03 &	0 &	3.02E-03 \\
$2\hskip-.6ex$ &
7.41E-03 &	3.78E-03 &	1.36E-03 &	1.51E-04 &	0 &	6.05E-04 \\
$3\hskip-.6ex$ &
5.19E-03 &	1.89E-03 &	4.08E-04 &	1.51E-05 &	0 &	1.21E-04 \\
$4\hskip-.6ex$ &
3.63E-03 &	9.46E-04 &	1.22E-04 &	1.51E-06 &	0 &	2.42E-05 \\
$5\hskip-.6ex$ &
2.54E-03 &	4.73E-04 &	3.67E-05 &	1.51E-07 &	0 &	4.84E-06 \\
$6\hskip-.6ex$ &
1.78E-03 &	2.36E-04 &	1.10E-05 &	1.51E-08 &	0 &	9.68E-07 \\
$7\hskip-.6ex$ &
1.24E-03 &	1.18E-04 &	3.31E-06 &	1.51E-09 &	0 &	1.93E-07 \\
$8\hskip-.6ex$ &
8.72E-04 &	5.91E-05 &	9.93E-07 &	1.51E-10 &	0 &	3.87E-08 \\
$9\hskip-.6ex$ &
6.10E-04 &	2.95E-05 &	2.97E-07 &	1.51E-11 &	0 &	7.75E-09 \\
$10\hskip-.6ex$ &
4.27E-04 &	1.47E-05 &	8.93E-08 &	1.51E-12 &	0 &	1.55E-09\\
\hline
\end{tabular}
\end{table}
\end{example}

\begin{example}\label{tridiag1}
We test Algorithm I (GEPP) on badly scaled tridiagonally matrix
$A$ generated by the \textsl{MATLAB} code
\medskip

\begin{verbatim}
randn('state',0)
n=10;m=5;
u=randn(n,1); v=randn(n-1,1);
A=diag(u)+diag(v,-1)+diag(v,1);
t=1e10;  A(m-1,m)=t;
end
\end{verbatim}

Random matrices of entries from the distribution $N(0,1)$. They
were generated by the \textsl{MATLAB} function "randn".
Before each usage the random number generator was reset to its
initial state.
Notice that only the element $A_{4,5}$ is very large (equals
$10^{10}$), hence the matrix $A$ is ill-conditioned.
The values of the componentwise stability error (\ref{co}) are
gathered in Table \ref{gepp_fo}.
Clearly the best results are obtained for $\omega=1$
(Wilkinson's original iterative refinement).
We don't display the forward error (\ref{fo}) and backward
stability error (\ref{ba}) because they were always small (of
order $\eps$).

\begin{table}
\caption{Values of the componentwise stability error (\ref{co})
for Algorithm I (GEPP) , where $A$ is the
$10\times 10$ tridiagonal matrix defined in Example
\ref{tridiag1} for $t= 10^{10}$. Here $\kappa(A)=7.74 \cdot
10^{10}$. }
\label{gepp_fo}
\begin{tabular}{lllllll}
\hline $\omega  \slash k$ & %
$0.3$ & $0.5$ & $0.7$ & %
$0.9$ & $1.0$ & %
$1.2$ \\
\hline 
$0\hskip-.6ex$ & 
1.02E-06 & 1.024E-06 & 1.02E-06 & 1.02E-06 & 1.02E-06 & 1.02E-06
\\
$1\hskip-.6ex$ &
7.15E-07 & 5.10E-07 & 3.06E-07 & 1.02E-07 & 1.15E-16 & 2.04E-07
\\
$2\hskip-.6ex$ &
5.00E-07 & 2.55E-07 & 9.19E-08 & 1.02E-08 & 1.15E-16 & 4.08E-08
\\
$3\hskip-.6ex$ &
3.50E-07 & 1.27E-07 & 2.75E-08 & 1.02E-09 & 1.15E-16 & 8.17E-09
\\
$4\hskip-.6ex$ &
2.45E-07 & 6.38E-08 & 8.27E-09 & 1.02E-10 & 1.15E-16 & 1.63E-09
\\
$5\hskip-.6ex$ &
1.71E-07 & 3.19E-08 & 2.48E-09 & 1.02E-11 & 1.15E-16 & 3.27E-10
\\
$6\hskip-.6ex$ &
1.20E-07 & 1.59E-08 & 7.44E-10 & 1.021E-12 & 1.15E-16 & 6.54E-11
\\
$7\hskip-.6ex$ &
8.41E-08 & 7.98E-09 & 2.23E-10 & 1.021E-13 & 1.15E-16 & 1.30E-11
\\
$8\hskip-.6ex$ &
5.89E-08 & 3.99E-09 & 6.70E-11 & 1.01E-14 & 1.15E-16 & 2.61E-12
\\
$9\hskip-.6ex$ &
4.12E-08 & 1.99E-09 & 2.01E-11 & 1.07E-15 & 1.15E-16 & 5.23E-13
\\
$10\hskip-.6ex$ &
2.88E-08 & 9.97E-10 & 6.034E-12 & 1.54E-16 & 1.15E-16 &
1.04E-13\\
\hline
\end{tabular}
\end{table}
\end{example}

\medskip

\begin{example}\label{BLUformacierz11hilb}
We generate a block matrix $A$ as in (\ref{block2}) using the
following \textsl{MATLAB} code.
\begin{verbatim}
m=8; n=2*m;
rand('state',0);
A=rand(n);
A(1:m,1:m)=hilb(m);
\end{verbatim}

The matrix $A$ is very well-conditioned, with the condition
number $\kappa(A)=2.0 8 \cdot 10^2$ but the block $(1,1)$ of $A$
is ill-conditioned: $\kappa(A_{11})= 4.75 \cdot 10^8$.
Here  $H=hilb(m)$ is a $m \times m$ Hilbert matrix defined by
\[
H=(h_{ij}), \,\,\, h_{ij}=\frac{1}{i+j-1},\,\, i, j = 1, \ldots,
m.
\]
 
The results are contained in Tables
\ref{macierz11hilb_fe}-\ref{macierz11hilb_co}.
 
\begin{table}
\caption{Values of the forward stability error (\ref{fo}) for
Algorithm II (BLU), where $A$ is the
$16 \times 16$ matrix defined in Example
\ref{BLUformacierz11hilb}. }
\label{macierz11hilb_fe}
\begin{tabular}{lllllll}
\hline $\omega  \slash k$ & %
$0.3$ & $0.5$ & $0.7$ & %
$0.9$ & $1.0$ & %
$1.2$ \\
\hline 
$0\hskip-.6ex$ & 
2.01E-10 & 2.01E-10 & 2.01E-10 & 2.01E-10 & 2.01E-10 & 2.01E-10
\\
$1\hskip-.6ex$ &
1.41E-10 & 1.00E-10 & 6.05E-11 & 2.01E-11 & 3.57E-17 & 4.03E-11
\\
$2\hskip-.6ex$ &
9.89E-11 & 5.04E-11 & 1.81E-11 & 2.01E-12 & 2.64E-17 & 8.07E-12
\\
$3\hskip-.6ex$ &
6.92E-11 & 2.52E-11 & 5.45E-12 & 2.01E-13 & 9.57E-18 & 1.61E-12
\\
$4\hskip-.6ex$ &
4.84E-11 & 1.26E-11 & 1.63E-12 & 2.01E-14 & 8.91E-18 & 3.23E-13
\\
$5\hskip-.6ex$ &
3.39E-11 & 6.31E-12 & 4.90E-13 & 2.01E-15 & 1.19E-17 & 6.46E-14
\\
$6\hskip-.6ex$ &
2.37E-11 & 3.15E-12 & 1.47E-13 & 1.96E-16 & 2.94E-17 & 1.29E-14
\\
$7\hskip-.6ex$ &
1.66E-11 & 1.57E-12 & 4.41E-14 & 2.62E-17 & 1.46E-17 & 2.58E-15
\\
$8\hskip-.6ex$ &
1.16E-11 & 7.88E-13 & 1.32E-14 & 3.71E-17 & 2.04E-17 & 5.12E-16
\\
$9\hskip-.6ex$ &
8.14E-12 & 3.94E-13 & 3.96E-15 & 5.36E-17 & 2.47E-17 & 9.83E-17
\\
$10\hskip-.6ex$ &
5.70E-12 & 1.97E-13 & 1.19E-15 & 2.70E-17 & 3.22E-17 &
4.84E-17\\
\hline
\end{tabular}
\end{table}

\medskip

\begin{table}
\caption{Values of the backward stability error (\ref{ba}) for
Algorithm II (BLU), where $A$ is the
$16 \times 16$ matrix defined in Example
\ref{BLUformacierz11hilb}. }
\label{macierz11hilb_ba}
\begin{tabular}{lllllll}
\hline $\omega  \slash k$ & %
$0.3$ & $0.5$ & $0.7$ & %
$0.9$ & $1.0$ & %
$1.2$ \\
\hline 
$0\hskip-.6ex$ & 
4.03E-09 & 4.03E-09 & 4.03E-09 & 4.03E-09 & 4.03E-09 & 4.03E-09
\\
$1\hskip-.6ex$ & 
2.82E-09 & 2.01E-09 & 1.21E-09 & 4.03E-10 & 1.90E-16 & 8.06E-10
\\
$2\hskip-.6ex$ & 
1.97E-09 & 1.00E-09 & 3.63E-10 & 4.03E-11 & 1.92E-16 & 1.61E-10
\\
$3\hskip-.6ex$ & 
1.38E-09 & 5.04E-10 & 1.08E-10 & 4.03E-12 & 1.43E-16 & 3.22E-11
\\
$4\hskip-.6ex$ & 
9.68E-10 & 2.52E-10 & 3.26E-11 & 4.03E-13 & 1.53E-16 & 6.45E-12
\\
$5\hskip-.6ex$ & 
6.78E-10 & 1.26E-10 & 9.80E-12 & 4.03E-14 & 1.44E-16 & 1.29E-12
\\
$6\hskip-.6ex$ & 
4.74E-10 & 6.30E-11 & 2.94E-12 & 4.01E-15 & 1.49E-16 & 2.58E-13
\\
$7\hskip-.6ex$ & 
3.32E-10 & 3.15E-11 & 8.82E-13 & 4.14E-16 & 1.63E-16 & 5.16E-14
\\
$8\hskip-.6ex$ & 
2.32E-10 & 1.575E-11 & 2.64E-13 & 1.14E-16 & 1.22E-16 & 1.03E-14
\\
$9\hskip-.6ex$ & 
1.62E-10 & 7.88E-12 & 7.93E-14 & 8.18E-17 & 1.66E-16 & 2.09E-15
\\
$10\hskip-.6ex$ & 
1.13E-10 & 3.94E-12 & 2.38E-14 & 1.44E-16 & 1.66E-16 &
5.16E-16\\
\hline
\end{tabular}
\end{table}

\begin{table}
\caption{Values of the componentwise backward stability error
(\ref{co}) for Algorithm II (BLU), where $A$ is the
$16 \times 16$ matrix defined in Example
\ref{BLUformacierz11hilb}. }
\label{macierz11hilb_co}
\begin{tabular}{lllllll}
\hline $\omega  \slash k$ & %
$0.3$ & $0.5$ & $0.7$ & %
$0.9$ & $1.0$ & %
$1.2$ \\
\hline 
$0\hskip-.6ex$ & 
7.88E-09 & 7.88E-09 & 7.88E-09 & 7.88E-09 & 7.88E-09 & 7.88E-09
\\
$1\hskip-.6ex$ &
5.51E-09 & 3.94E-09 & 2.36E-09 & 7.88E-10 & 4.19E-16 & 1.57E-09
\\
$2\hskip-.6ex$ &
3.86E-09 & 1.97E-09 & 7.09E-10 & 7.88E-11 & 4.61E-16 & 3.15E-10
\\
$3\hskip-.6ex$ &
2.70E-09 & 9.85E-10 & 2.12E-10 & 7.88E-12 & 3.07E-16 & 6.30E-11
\\
$4\hskip-.6ex$ &
1.89E-09 & 4.92E-10 & 6.38E-11 & 7.88E-13 & 3.07E-16 & 1.26E-11
\\
$5\hskip-.6ex$ &
1.32E-09 & 2.46E-10 & 1.91E-11 & 7.89E-14 & 2.79E-16 & 2.52E-12
\\
$6\hskip-.6ex$ &
9.27E-10 & 1.23E-10 & 5.74E-12 & 7.87E-15 & 2.27E-16 & 5.04E-13
\\
$7\hskip-.6ex$ &
6.49E-10 & 6.15E-11 & 1.72E-12 & 7.95E-16 & 3.07E-16 & 1.01E-13
\\
$8\hskip-.6ex$ &
4.54E-10 & 3.07E-11 & 5.17E-13 & 2.25E-16 & 2.21E-16 & 2.03E-14
\\
$9\hskip-.6ex$ &
3.18E-10 & 1.53E-11 & 1.54E-13 & 1.89E-16 & 3.04E-16 & 4.15E-15
\\
$10\hskip-.6ex$ &
2.22E-10 & 7.69E-12 & 4.65E-14 & 3.78E-16 & 3.41E-16 & 1.02E-15
\\
\hline
\end{tabular}
\end{table}
\end{example}

\medskip

Based on the numerical results of this section, we conclude that
one step of Wilkinson's iterative refinement ($\omega=1$) is
usually be enough to yield small errors (\ref{fo})--(\ref{co}).
However, iterative refinement with a relaxation $\omega$ which
is not close to $1$, can require much more steps than
Wilkinson's iterative refinement.
Therefore, the choice $\omega=1$ is the best choice from the
point of numerical stability.


\medskip

\begin{thebibliography}{0}
\bibitem{arioli} M.~Arioli, J.~Demmel, I.S.~Duff,
Solving sparse linear systems with sparse backward error,
SIAM J. Matrix Anal. Appl. 10 (1989), 165--190.

\bibitem{Arioli} M. Arioli, J.Scott,   Chebyshev acceleration of iterative refinement, Numer. Algor.  66 (2014), 591–-608.

\bibitem{bjorck1} {\AA}. Bj\"{o}rck, Iterative refinement and
realiable computing. In Reliable Numerical Computation,
M.~G.~Cox and
S.~J.~Hammarling, editors, Oxford University Press  (1990),
249--266.

\bibitem{bjorck2} {\AA}.~ Bj\"{o}rck, Iterative refinement of linear least squares solutions I,
BIT 7 (1967),  257--278.

\bibitem{bjorck3} {\AA}.~ Bj\"{o}rck, Iterative refinement of linear least squares solutions II,
BIT 8 (1968),  8--30.


\bibitem{Demmel: 95} J.W.Demmel, N.J.Higham, and R.S.Schreiber,
Stability of block LU factorization,
Numer. Linear Algebra Appl. 12 (1995), 173--190.

\bibitem{Foster} L.~V.~Foster, Gaussian elimination with partial
pivoting can fail in practice, SIAM J.Matrix Anal. Appl. 15 (4) (1994), 1354--1362.

\bibitem{higham1} N.J.Higham, Iterative refinement enhances the
stability of QR factorization methods
for solving linear equations, BIT 31 (1991),  447--468.

\bibitem{higham2} N.J.Higham,  Iterative refinement for linear systems and LAPACK,   IMA J. Numer. Anal., 17(4) (1997),  495--509.

\bibitem{higham:96} N.J.Higham, {\it Accuracy and Stability of
Numerical
Algorithms}, Second Edition, SIAM, Philadelphia (2002).

\bibitem{jankowski} M.Jankowski, H.Wo\'{z}niakowski, Iterative
refinement implies numerical stability,
BIT 17 (1977), 303--311.

\bibitem{Cleve} C.B.Moler, Iterative refinement in floating
point, J.Assoc. Comput. Mach. 14 (2) (1967), 316--321.

\bibitem{Miro} M.Rozlo\v{z}n\'{\i}k, A.Smoktunowicz and J.Kopal, A note on iterative refinement fo rseminormal equations, 
Applied Numerical Mathematics 75 (2014), 167--174.

\bibitem{skeel} R.D.Skeel, Iterative refinement implies
numerical stability for Gaussian elimination,
Math. Comp. 35 (1980),  817--832.

\bibitem{smok3} Alicja Smoktunowicz and Agata Smoktunowicz,
Iterative refinement techniques for solving block linear systems
of equations, Applied Numerical Mathematics 67 (2013), 220-229.

\bibitem{Wilkinson} J.~H.~Wilkinson, {\it The Algebraic
Eigenvalue Problem}, Oxford University Press (1965).

\bibitem{Wu} X.Wu and Z.Wang, A new iterative refinement with
roundoff error analysis, Numer. Linear Algebra Appl. 18 (2011), 275--282.

\bibitem{Zlatev} Z. Zlatev, Use of iterative refinement in the
solution of sparse linear systems, SIAM J. Numer. Anal., 19(2)
(1982), 381–-399.
\end{thebibliography}
\end{document}